\documentclass[12pt]{amsart}
\usepackage[top=3cm, bottom=2cm, left=4.0cm, right=3.5cm]{geometry}

\theoremstyle{definition}

\usepackage{epsfig,color}

\makeatother

\theoremstyle{definition}

\def\fnum{equation} 
\newtheorem{Thm}[\fnum]{Theorem}

\numberwithin{equation}{section}

\usepackage[hidelinks]{hyperref}

\begin{document}
\title
[Rigidity for the logarithmic Sobolev inequality]
{Rigidity for the logarithmic Sobolev inequality on complete metric measure spaces}

\author[Franciele Conrado]{Franciele Conrado} \address{Instituto de Matemática e Estatística \\ 
Universidade Federal Fluminense\\ 
Niterói-RJ, Brazil}\email{franciconradomat@gmail.com} \thanks{The author was Partially supported by CNPq-Brazil} 

\begin{abstract}
In this work, we study the rigidity problem for the logarithmic Sobolev inequality on a complete metric measure space $(M^n,g,f)$ with Bakry–Émery Ricci curvature satisfying $Ric_f\geq \frac{a}{2}g$, for some $a>0$. We prove that if equality holds then $M$ is isometric to $\Sigma\times \mathbb{R}$ for some complete $(n-1)$-dimensional Riemannian manifold $\Sigma$ and by passing an isometry, $(M^n,g,f)$ must split off the Gaussian shrinking soliton $(\mathbb{R}, dt^2, \frac{a}{2}|.|^2)$. This was proved in 2019 by Ohta and Takatsu in \cite{ohta}. In this paper, we prove this  rigidity result using a different method.
\end{abstract}

\maketitle

\section{Introduction}\label{intro}

\medskip

Let $(M^n,g)$, $n\geq 2$, be a compact $n$-dimensional Riemannian manifold. It is well-known that the spectrum of the Laplacian $\Delta$ is discrete. The relation between the spectrum of the Laplacian and geometric quantities has been a topic of continued interest. One of such relation was given by Lichnerowicz in \cite{L}, which proved that if the Ricci curvature satisfies $Ric\geq (n-1)ag$, where $a$ is a positive constant, the first nonzero eigenvalue of the Laplacian $\Delta$ satisfies $\lambda_1(\Delta)\geq na$. Recall that for the $n$-dimensional round sphere of radius $\frac{1}{\sqrt{a}}$ we have that $Ric=(n-1)ag$ and $\lambda_1(\Delta)=na$. Subsequently, Obata in \cite{obata} proved that $\lambda_1(\Delta)= na$ if and only if $M$ is isometric to an $n$-dimensional round sphere of radius $\frac{1}{\sqrt{a}}$. We have that the Poincaré inequality
\begin{equation}\label{ee}\int_M\rho^2 dv \leq \frac{1}{na}\int_M |\nabla \rho|^2 dv\end{equation}
\noindent holds for every locally Lipschitz function $\rho\in L^2(M)$ with $\int_M \rho dv=0$. Here, $dv$ denotes the volume element of $(M^n,g)$ and $L^2(M)$ denotes the space of square-integrable functions on $M$. Observe that $\lambda_1(\Delta)=na$ if and only if the equality holds in the Poincaré inequality (\ref{ee}).

\medskip

Let $(M^n,g,f)$ be a complete smooth metric measure space, i.e., a complete $n$-dimensional Riemannian manifold $(M^n,g)$ together with a function $f\in C^{\infty}(M)$. Denote by $\mu$ the measure induced by the weighted volume element $e^{-f}dv$. We define the Bakry–Émery Ricci curvature by $Ric_f= Ric+\nabla^2f$ and the drifted Laplacian by $\Delta_f=\Delta-\langle \nabla f, \nabla \cdot \rangle$, where $Ric$ is the Ricci curvature tensor, $\nabla^2 f$ is the Hessian of the function $f$ and $\Delta$ is the Laplacian on $(M^n,g)$.

\medskip

Assume that $Ric_f\geq \frac{a}{2}g$, for some $a>0$. In this case, Morgan proved in \cite{morgan} that the weighted volume $\mu(M)=\int_Me^{-f}dv$ is finite. Consequently, by results proved in \cite{BE}, the logarithmic Sobolev inequality
\begin{equation}\label{e1}\int_M\rho \log \rho d\mu\leq \frac{1}{a}\int_M\frac{|\nabla \rho|^2}{\rho} d\mu \end{equation}
\noindent holds for every nonnegative locally Lipschitz function $\rho:M\rightarrow \mathbb{R}$ with $\int_M\rho d\mu=\mu(M)$. In particular, the Poincaré inequality
\begin{equation}\label{e}\int_M\rho^2 d\mu \leq \frac{2}{a}\int_M |\nabla \rho|^2 d\mu \end{equation}
\noindent holds for every locally Lipschitz function $\rho\in L^2(M,\mu)$ with $\int_M \rho d\mu=0$. Here, $L^2(M,\mu)$ denote the space of square-integrable functions on $M$ with respect to the measure $\mu$.  It is well-known that the logarithmic Sobolev inequality (\ref{e1}) and the finiteness of $\mu(M)$ imply the discreteness of the spectrum of $\Delta_f$.

\medskip

In 2015, Cheng and Zhou proved in \cite{XCDZ} an analogous result to the Lichnerowicz–Obata’s theorem for a complete smooth metric measure space $(M^n,g,f)$ with Bakry–Émery Ricci curvature satisfying $Ric_f\geq \frac{a}{2}g$, for some $a>0$. They proved that the first nonzero eigenvalue of the drifted Laplacian $\Delta_f$ satisfies $\lambda_1(\Delta_f)\geq \frac{a}{2}$ and that equality holds if and only if $M$ is isometric to $\Sigma\times \mathbb{R}$ for some complete $(n-1)$-dimensional Riemannian manifold $\Sigma$ and by passing an isometry, $(M^n,g,f)$ must split off the Gaussian shrinking soliton $(\mathbb{R}, dt^2, \frac{a}{2}|.|^2)$.

\medskip

Observe that $\lambda_1(\Delta_f)=\frac{a}{2}$ if and only if the equality holds in the Poincaré inequality (\ref{e}). Furthermore, since the logarithmic Sobolev inequality implies the Poincaré inequality, then the equality in (\ref{e1}) is
a weaker assumption than the equality in (\ref{e}). Hence, studying the rigidity problem for the logarithmic Sobolev inequality is a very interesting problem. 

\medskip

Recently, Ohta and Takatsu in \cite{ohta} studied the rigidity for the logarithmic Sobolev inequality in complete smooth metric measure space $(M^n,g,f)$ with Bakry–Émery Ricci curvature satisfying $Ric_f\geq \frac{a}{2}g$, for some $a>0$. They used the needle decomposition method introduced on Riemannian manifolds by Klartag to prove that if the equality holds in (\ref{e1}), then the same rigidity result obtained by Cheng and Zhou holds, i.e., $M$ is isometric to $\Sigma\times \mathbb{R}$ for some complete $(n-1)$-dimensional Riemannian manifold $\Sigma$ and by passing an isometry, $(M^n,g,f)$ must split off the Gaussian shrinking soliton $(\mathbb{R}, dt^2, \frac{a}{2}|.|^2)$.  In this paper, we prove this rigidity using a different method, namely, we use the ideas of Cheng and Zhou presented in \cite{XCDZ}. More precisely, we prove the following:

\begin{Thm}\label{tt}
Let $(M^n,g,f)$ be a complete smooth metric measure space such that $Ric_f\geq \frac{a}{2}g$, for some $a>0$, and the measure $\mu$ induced by the weighted volume element $e^{-f}dv$ is a probability measure. Assume 
$$\int_M\rho \log \rho d\mu = \frac{1}{a}\int_M\frac{|\nabla \rho|^2}{\rho} d\mu < +\infty $$
\noindent for some nonconstant and locally Lipschitz function $\rho:M\rightarrow [0,+\infty)$ with $\int_M\rho d\mu=1$. Then
\begin{itemize}
\item $M$ is noncompact;
\item $M$ is isometric to $\Sigma\times \mathbb{R}$, where $(\Sigma,g_{\Sigma})$ is an $(n-1)$-dimensional complete Riemannian manifold;
\item By passing an isometry, for $(x,t)\in \Sigma\times \mathbb{R}$, we have that
$$f(x,t)=\frac{a}{4}t^2+f(x,0).$$
Moreover, $Ric_{\tilde{f}}^{\Sigma}\geq \displaystyle\frac{a}{2}g_{\Sigma}$, where $\tilde{f}=f(.,0)$.
\end{itemize}
\end{Thm}

\medskip


\subsection*{Acknowledgments} The author would like to thank Professor Detang Zhou for his interest and helpful discussions.

\section{Proof of the theorem \ref{tt}}

\begin{proof}
Define the function $u:=\sqrt{\rho}$. Note that $u\in H^1(M,\mu)$, where $H^1(M,\mu)$ denotes the space of functions in $L^2(M,\mu)$ whose gradient is square-integrable with respect to the measure $\mu$. Furthermore, note that $\int_M u^2 d\mu=1$ and
\begin{equation}\label{e2}\int_M u^2 \log u^2 d\mu = \frac{4}{a}\int_M |\nabla u|^2 d\mu.\end{equation}
 
Let $F,G:H^1(M,\mu)\rightarrow  \mathbb{R}$ be functionals defined by
$$F(w)=\int_M \left(\frac{4}{a} |\nabla w|^2-w^2\log w^2\right) d\mu \ \ \ \text{and} \ \ \ G(w)=\int_M w^2 d\mu.$$

It follows  from the Logarithmic Sobolev Inequality (\ref{e1}) and the equality (\ref{e2}) that $F(w)\geq 0=F(u)$ for every $w\in G^{-1}(1)$. By the method of Lagrange multipliers, there is $\lambda\in \mathbb{R}$ such that $DF_u=\lambda DG_u$. Consequently,
$$ \int_Mw\left(-\frac{8}{a}\Delta_fu-4u\log u-2u\right)d\mu =2\lambda \int_Mwu d\mu,$$ 
\noindent for every $w\in H^1(M,\mu)$. It follows that
$$\int_Mw\left(\Delta_fu+\frac{a}{4}u(2\log u+1+\lambda)\right)d\mu =0,$$ 
\noindent for every $w\in H^1(M,\mu)$. This implies that $u$ is a weak solution to
\begin{equation}\label{e3}
\Delta_fu+\frac{a}{4}u(\log u^2+1+\lambda)=0.
\end{equation}

Since $F(u)=0$  and $u$ is a weak solution to (\ref{e3}), we obtain $\lambda=-1$. Consequently, we have that $u$ is a weak solution to
\begin{equation}\label{e4}\Delta_fu+\frac{a}{4}u\log u^2=0.\end{equation}

Using results from regularity theory we obtain that $u\in C^{\infty}(M)$ and then $u$ is a classical solution to (\ref{e4}). Define $v:=\log u$ and note that
\begin{eqnarray*}
\Delta_fu+\frac{a}{4}u\log u^2 &=& \Delta_f(e^v)+\frac{a}{4}e^v\log (e^{2v})\\
&=&\Delta (e^v)-\langle \nabla f, \nabla (e^v)\rangle +\frac{a}{2}ve^v\\
&=& e^v\Delta v + e^v |\nabla v|^2-e^v\langle \nabla f, \nabla v \rangle +\frac{a}{2}ve^v\\
&=& e^v\left(  \Delta_f v + |\nabla v|^2+\frac{a}{2}v\right).\\
\end{eqnarray*}

It follows that
\begin{equation}\label{e5}
 \Delta_f v + |\nabla v|^2+\frac{a}{2}v=0.
\end{equation}

By $Ric_f\geq \displaystyle\frac{a}{2}g$, (\ref{e5}) and  the weighted Bochner formula:
$$\frac{1}{2}\Delta_f(|\nabla v|^2)=|\nabla^2v|^2+\langle \nabla v, \nabla (\Delta_f v)\rangle+ Ric_f(\nabla v,\nabla v)$$
\noindent we have that
\begin{equation}\label{e6}\frac{1}{2}\Delta_f(|\nabla v|^2)\geq |\nabla^2v|^2-\langle \nabla v, \nabla (|\nabla v|^2)\rangle .\end{equation}

Note that
\begin{equation}\label{e7}\Delta_f(|\nabla v|^2)+ 2\langle \nabla v, \nabla (|\nabla v|^2)\rangle = \Delta_h (|\nabla v|^2),\end{equation}
\noindent where $h:=f-2v$. By (\ref{e6}) and (\ref{e7}) we obtain
\begin{equation}\label{e8}\Delta_h (|\nabla v|^2) \geq 2|\nabla^2v|^2.\end{equation}

Denote by $\nu$ the measure induced by the weighted volume element  $e^{-h}dv$. Note that $\int_M|\nabla v|^2d\nu < +\infty$ because $u\in H^1(M,\mu)$ and
$$\int_M |\nabla v|^2d\nu= \int_M |\nabla u|^2d\mu.$$

Fix a point $p\in M$. For each $k\in \mathbb{N}$, consider $B_k$ the metric ball in $M$ of radius $k$ centered at $p$. Consider a nonnegative function $\phi_k\in C^{\infty}(M)$ such that $\phi_k\equiv 1$ on $B_k$, $|\nabla \phi_k|\leq 1$ on $B_{k+1}\setminus B_k$ and $\phi_k\equiv 0$ on $M\setminus B_{k+1}$. By (\ref{e8}), we have that
\begin{equation}\label{e9}\phi_k^2\Delta_h (|\nabla v|^2) \geq 2\phi_k^2|\nabla^2v|^2.\end{equation}

Note that
\begin{eqnarray*}
\left|\int_M \phi_k^2\Delta_h (|\nabla v|^2) d\nu \right| & = & \left|-2\int_M \phi_k \langle \nabla \phi_k,\nabla (|\nabla v|^2)\rangle d\nu \right|\\
& = & \left|2\int_M \phi_k \nabla \phi_k(|\nabla v|^2) d\nu\right|\\
& = & \left|2\int_M 2\phi_k \langle \nabla_{\nabla \phi_k}\nabla v, \nabla v\rangle d\nu\right|\\
& = & \left|2\int_M 2\phi_k \nabla^2v(\nabla \phi_k,\nabla v) d\nu\right|\\
& \leq & 2\int_M |2\phi_k \nabla^2v(\nabla \phi_k,\nabla v)| d\nu.\\
\end{eqnarray*}

Let $\{e_1,\dots ,e_n\}$ be a local orthonormal frame on $(M^n,g)$. Fix a number $0<\epsilon<1$. By Young’s inequality $2ab\leq \epsilon a^2+\epsilon^{-1}b^2$, we have
\begin{eqnarray*}
|2\phi_k \nabla^2v(\nabla \phi_k,\nabla v)| &\leq & \sum_{i,j=1}^n 2\phi_k |e_i(\phi_k)| |e_j(v)| |(\nabla^2v)_{ij}| \\
&\leq & \sum_{i,j=1}^n \left[\epsilon \phi_k^2 (\nabla^2v)^2_{ij} + \frac{1}{\epsilon}(e_i(\phi_k))^2(e_j(v))^2\right]\\
&=&  \epsilon\phi_k^2|\nabla^2v|^2 +\frac{1}{\epsilon}|\nabla \phi_k|^2 |\nabla v|^2.\\
\end{eqnarray*}

It follows that 
$$\left|\int_M \phi_k^2\Delta_h (|\nabla v|^2) d\nu \right| \leq 2\epsilon\int_M \phi_k^2|\nabla^2v|^2 d\nu +\frac{2}{\epsilon}\int_M |\nabla \phi_k|^2 |\nabla v|^2d\nu .$$

So (\ref{e9}) implies that
$$(1-\epsilon)\int_M \phi_k^2|\nabla^2v|^2 d\nu \leq \frac{1}{\epsilon}\int_M |\nabla \phi_k|^2 |\nabla v|^2 d\nu .$$

Consequently,
\begin{eqnarray*}
\int_{B_k} |\nabla^2v|^2 d\nu &=& \int_{B_k} \phi_k^2|\nabla^2v|^2 d\nu\\
& \leq & \int_M \phi_k^2|\nabla^2v|^2 d\nu \\
&\leq & \frac{1}{\epsilon(1-\epsilon)}\int_M |\nabla \phi_k|^2 |\nabla v|^2 d\nu\\
& = & \frac{1}{\epsilon(1-\epsilon)}\int_{B_{k+1}\setminus B_k} |\nabla \phi_k|^2 |\nabla v|^2 d\nu\\
&\leq & \frac{1}{\epsilon(1-\epsilon)}\int_{B_{k+1}\setminus B_k} |\nabla v|^2 d\nu.\\
\end{eqnarray*}

This implies that
\begin{equation}\label{e11} \int_{B_k} |\nabla^2v|^2 d\nu \leq  \frac{1}{\epsilon(1-\epsilon)}\int_{B_{k+1}\setminus B_k} |\nabla v|^2 d\nu .\end{equation}

For each Borel set $A\subset M$ we denote by $\chi_A$ the characteristic function of $A$. Since $\int_M|\nabla v|^2d\nu < +\infty$, $|\nabla v|^2\chi_{B_{k+1}\setminus B_k} \rightarrow 0$ and $|\nabla v|^2\chi_{B_{k+1}\setminus B_k}\leq |\nabla v|^2$ for every $k\in\mathbb{N}$, by the dominated convergence theorem we have
\begin{equation}\label{e12}\lim_{k \to +\infty} \int_{B_{k+1}\setminus B_k} |\nabla v|^2 d\nu =0.\end{equation}

Furthermore, since $|\nabla^2v|^2\chi_{B_k}\rightarrow |\nabla^2v|^2$ e $|\nabla^2v|^2\chi_{B_k}\leq |\nabla^2v|^2\chi_{B_{k+1}}$ for every $k\in\mathbb{N}$, by the monotone convergence theorem we have
\begin{equation}\label{e13} \lim_{k \to +\infty} \int_{B_k} |\nabla^2v|^2 d\nu = \int_M |\nabla^2v|^2 d\nu .\end{equation}

From (\ref{e11}), (\ref{e12}) and (\ref{e13}), we obtain
$$\int_M |\nabla^2v|^2 d\nu\leq 0.$$

Hence $\nabla^2v = 0$. This implies that $M$ is noncompact, because $v$ is a nonconstant harmonic function on $M$. Furthermore, we have that $\nabla v$ is a nontrivial parallel vector field on $M$. It follows that $(M^n,g)$ is isometric to $(\Sigma\times \mathbb{R}, g_{\Sigma}+dt^2)$, where $(\Sigma,g_{\Sigma})$ is an $(n-1)$-dimensional complete Riemannian manifold. By passing an isometry, we may assume that  $(M^n,g)=(\Sigma\times \mathbb{R}, g_{\Sigma}+dt^2)$ and $v$ restrict to $\Sigma\times \{t\}$ is constant, for every $t\in \mathbb{R}$. Consider $(x_1, \cdots, x_{n-1},t)$ local coordinates on $\Sigma\times \mathbb{R}=M$ such that $(x_1, \cdots, x_{n-1})$ are local coordinates on $\Sigma$. Note that $\nabla v=\lambda\partial_t$, for some nonzero constant $\lambda\in\mathbb{R}$, because $\nabla v$ is a nontrivial parallel vector field and $\partial_{x_i}(v)=0$ for every $i=1,\cdots , n-1$. Without loss of generality, suppose that $\Sigma=v^{-1}\left( -\frac{2\lambda^2}{a}\right)$. Then, 
$$v=\lambda t-\frac{2\lambda^2}{a}.$$

Since $\Delta v=0$, it follows from (\ref{e5}) that
$$-\langle \nabla v,\nabla f\rangle + |\nabla v|^2+\frac{a}{2}v=0.$$

This implies that
$$-\lambda\partial_t(f)+\lambda^2+\frac{a}{2}\left(\lambda t-\frac{2\lambda^2}{a}\right)=0.$$

Consequently, $\partial_t(f)=\displaystyle\frac{a}{2}t$. Then
$$f(x,t)=\frac{a}{4}t^2+f(x,0)$$
\noindent for every $(t,x)\in\Sigma\times \mathbb{R}=M.$ Furthermore, denote by $\tilde{f}=f(.,0)$ and note that, if $X$ is a smooth vector field on $\Sigma$, we have 
$$Ric^{\Sigma}_{\tilde{f}}(X,X)=Ric_f(X,X)\geq \frac{a}{2}|X|^2_g =\frac{a}{2}|X|^2_{g_{\Sigma}}.$$

Hence,
$$Ric^{\Sigma}_{\tilde{f}}\geq \frac{a}{2} g_{\Sigma}.$$
\end{proof}

\bibliographystyle{abbrv}
\bibliography{bib}

\end{document}